\documentclass{article}

\usepackage{amssymb}
\usepackage{amsmath}
\usepackage{mathtools}

\usepackage[english]{babel}

\newtheorem{theorem}{Theorem}[section]
\newtheorem{lemma}[theorem]{Lemma}

\newtheorem{remark}{\it Remark\/}

\newcommand{\fiz}{ { \textstyle{ \frac{1}{2} } } }

\newcommand{\II}{{\mathcal I}}
\newcommand{\JJ}{{\mathcal J}}
\newcommand{\KK}{{\mathcal K}}
\newcommand{\bH}{{\bf H}}
\newcommand{\bS}{{\bf S}}

\newcommand{\PP}{{\bf P}}

\begin{document}
\title{A Quaternionic Structure as a Landmark for Symplectic Maps}

\author{Hugo Jim\'enez-P\'erez 
\thanks{Funding: The author was supported by the \emph{Fondation du Coll\`ege de France},
\emph{Total} (RCN PU14150472), and the ERC Advanced Grant 
WAVETOMO (RCN 99285). 
The author declares that he has no conflict of interest.
}}
\maketitle

\begin{abstract}
We use a quaternionic structure on the product of two symplectic manifolds for 
relating Liouvillian forms with linear symplectic maps obtained by the symplectic Cayley's transformation.

\end{abstract}
\section{Introduction}
\label{}
One of the main difficulties for constructing symplectic maps by the method of
generating functions is the resolution of the Hamilton-Jacobi equation.
Instead of solving such an equation,
we consider a local quaternionic structure 
on the symplectic product manifold and the three different symplectic 
forms induced by this structure. Symplectic maps are constructed using the 
primitive Liouvillian forms related to those symplectic forms.

\section{Local symplectic maps from Liouvillian forms}
\label{sec:maps}
Let $(M,\omega)$ be a $2n$-dimensional symplectic manifold with symplectic form $\omega$.
A \emph{symplectomorphism} is a diffeomorphism $\phi:(M,\omega)\to (M,\omega)$
preserving the symplectic structure $\phi^*\omega= \omega$, where the star 
stands for the  \emph{pull-back} of differential forms. 
When the symplectic structure has a global primitive linear form $\theta$, 
then $(M,d\theta)$ is called an \emph{exact symplectic manifold}. Main representatives are cotangent bundles 
$(T^*\mathcal Q,d\alpha)$ which possesses a \emph{canonical} or \emph{tautological form} $\alpha$ called the \emph{Liouville form}. 
We define a \emph{Liouvillian form} in an exact symplectic manifold, as 
a primitive form $\theta\in\Omega^1(M)$ for the symplectic structure $\omega=d\theta$.
A \emph{Liouville vector field} $Z$ is given as the symplectic dual of $\theta$ 
by the implicit equation $\theta = (i_Z\omega)$.

We are interested in symplectic maps for constructing symplectic integrators, 
then we consider maps close to the identity 
defined on convex balls $\mathcal B\subset (M,\omega)$
containing at the same time, the source and the target points $z_0,z_\tau\in \mathcal B$.
We consider as hypothesis that these points belong to the flow $\phi^t=\phi(t)$ of 
a symplectic vector field, and by the Poicar\'e's lemma to a local Hamiltonian 
vector field. It implies that the 
segment of the Hamiltonian flow connecting them $z_\tau=\phi^\tau( z_0)$ 
is an embedded segment of curve $\phi^t([0,\tau])\hookrightarrow \mathcal B$. 
In a convex ball, there always exist primitive 1-forms $\theta$
by Poincar\'e's lemma and consequently we can apply this procedure locally on any symplectic manifold. 

\subsection {The product manifold}
As usual, define the product $\PP= M_1\times M_2$ of two 
copies of an exact symplectic manifold $(M,\omega=d\theta)$, which we denote by
$(M_1,\omega_1)$ and $(M_2,\omega_2)$, respectively. Each copy corresponds to the flow 
of a (Hamiltonian) system at two different times $t=0$ and $t=\tau$ for small $\tau$.
The canonical projections 
$\pi_i:\PP\to M_i$ for $i=1,2$ let us define the forms $\theta_\ominus$ and $\omega_{\ominus}$ on $\PP$ by,
$\theta_{\ominus} = \pi_1^*\theta_1 - \pi_2^*\theta_2,$ and $\omega_{\ominus} = \pi_1^*\omega_1 - \pi_2^*\omega_2$.
It is well known that $(\PP,\omega_{\ominus})$ is a symplectic manifold of dimension $4n$ \cite{LM87}. 
The graph of any symplectic map $\phi:(M_1,\omega_1)\to (M_2,\omega_2)$, defined by 
$$\Gamma_\phi = \left\{ \left( z,\phi(z) \right)\in \PP\ |\ z\in M_1,\phi(z)\in M_2 \right\},$$
is a Lagrangian submanifold in $(\PP,\omega_\ominus)$. 

An embedding of a half-dimensional manifold $\jmath:N\hookrightarrow (M,\omega)$
into a symplectic manifold is called \emph{Lagrangian} if 
$\jmath^*\omega\equiv 0$. 
Consequently, we can consider $\Gamma_\phi$ as an embedding $\jmath:\Lambda\hookrightarrow \PP$ with 
$\jmath(\Lambda)=\Gamma_\phi$ satisfying $\jmath^*\omega_\ominus = 0$. 
In addition to the symplectic form $\omega_\ominus$, for every $x\in \PP$ there exists an induced endomorphism 
on $T_x\PP$ which becomes the  complex structure associated to $\omega_\ominus$ given by 
$ J_\ominus = J_1 \oplus J_2^T$, where $J_i$, are the complex structures
associated to $\omega_i$, $i=1,2$.

The linear form $\alpha=\jmath^*\theta_\ominus$ 
on $\Lambda$ is closed since its differential satisfies  
$$d\alpha =\jmath^*d\theta_\ominus= \jmath^*\omega_\ominus \equiv 0.$$ 
Applying Poincare's lemma, $\alpha$ is locally exact on  $\Lambda$ and 
there (locally) exists a  function $S:\Lambda\to\mathbb R$ defined on
$\Lambda$ such that its differential concides with the pullback of $\theta_{\ominus}$ to $\Lambda$,
\emph{i.e.}
$dS  = \alpha = \jmath^*\theta_\ominus$. 
The function $S:\Lambda\to\mathbb R$ is called a \emph{generating 
function} for the symplectic map $\phi:(M,\omega)\to(M,\omega)$. In fact the generating 
function is a function $\hat S:\PP\to\mathbb R$ defined on the image 
$\jmath(\Lambda)\subset\PP$ and the function $S$ is 
the composition $S\equiv\hat S\circ \jmath :\Lambda\to\mathbb R$.
A related generating function 
$F:T^*(Q_1\times Q_2)\to \mathbb R$ is given by the pullback of $S$ under the 
diffeomorphism  $\epsilon:\PP\to T^*(Q_1\times Q_2)$, such that
$dF=\epsilon^*(dS)$.

Given a Lagrangian embedding $\jmath:\Lambda\hookrightarrow (\PP,\omega_\ominus)$,
there exists an open neighborhood $\Lambda\subset U\subset \PP$ around $\Lambda$
and a projection $\pi : U\to \Lambda$, such that the composition 
$\Lambda \xhookrightarrow{\jmath} U \stackrel{\pi}{\to} \Lambda$
satisfies $\pi\circ\jmath=id_\Lambda$. This fact is just  
Weinstein's theorem saying that $U$ is locally symplectomorphic to an 
open neighborhood of the zero section in  $T^*\Lambda$ \cite{Wei71}. 
A Liouvillian form $\theta$ on $(\PP,\omega_\ominus=d\theta)$ is related to the generating 
function $S:\Lambda\to \mathbb R$ 
by the identity $dS=\jmath^*\theta$, and it satisfies $\ker \theta \subset \ker \pi^*(dS)$, equivalently $ \ker\theta \subset \jmath_*(T\Lambda)$.
The last relation is all we need to know to construct symplectic maps from Liouvillian forms.

\subsection{A quaternionic structure on the product manifold}
The method of generating functions uses two different 
symplectic structures on $\PP$, usually denoted by $\omega_\ominus$ and $\omega_\oplus$,
for working with Lagrangian submanifolds
\cite{MS17,DaS03}. 
It implicitly uses a \emph{twist} diffeomorphism considered the canonical isomorphism 
for cotangent bundles relating 
$T^*\mathcal Q_1\times T^*\mathcal Q_2\cong T^*(\mathcal Q_1\times\mathcal Q_2).$

For the construction of symplectic maps, a different twist diffeomorphism is applied
solving an alternative Hamiltonian system \cite{DaS03}. This diffeomorphism relates the 
product manifold with the double cotangent bundle
$T^*\mathcal Q\times T^*\mathcal Q\cong T^*(T^*\mathcal Q),$
and defines a projection by composition 
$$T^*\mathcal Q_1\times T^*\mathcal Q_2  \stackrel{\Phi}{\to} T^*(T^*\mathcal Q) 
\stackrel{\pi_{T^*\mathcal Q}}{\to} T^*\mathcal Q.$$
The way we select the twist $\Phi$ will define a different projection 
which, by the way, it determines a particular type of generating function. 

In this paper, we avoid the twist diffeomorphisms and the uncomfortable situation of working 
with different symplectomorphic manifolds. Instead, we consider only
the product manifold $\PP$, and we define an  
\emph{almost quaternionic} or \emph{almost hypercomplex structure} on $\PP$ given by
$\{(I_{4n},\II,\JJ,\KK\}\subset End(T\PP)$ \cite{Boy88}, which
induces the local geometry of $(\PP,\omega_\ominus)$, 
$(T^*(Q_1\times Q_2),\omega_\oplus)$ and $(T^*M,\omega_{can})$.  In Darboux coordinates, we have 
the matricial representation\footnote{We use $\II=J_{4n}^T$ in accordance to complex geometry. 
See the discussion in \cite[Rmk. 3.1.6]{MS17}.} 
\begin{eqnarray}
   \II = 
   \left(\begin{smallmatrix}
        0_{2n} & -I_{2n}\\
        I_{2n} & 0_{2n}
   \end{smallmatrix}\right), \quad 
    \JJ = 
   \left(\begin{smallmatrix}
        J_{2n} & 0_{2n}  \\
        0_{2n} & J_{2n}^T
   \end{smallmatrix}\right) 
    \qquad{\rm and }\qquad
    \KK =    
   \left(\begin{smallmatrix}
        0_{2n} & J_{2n}\\
        J_{2n} & 0_{2n}
   \end{smallmatrix}\right),
\end{eqnarray}
satisfying
\begin{eqnarray}
   \II^2 = \JJ^2 = \KK^2  = \II\JJ\KK = -I_{4n},\quad  \II\JJ= \KK, 
   \quad \JJ\KK=\II, \quad \KK\II= \JJ.
\end{eqnarray}
We obtain an equivalent framework to the usual one, and it is easy to prove that it just 
corresponds to a relabeling of coordinates. 

Let $g$ be the 
Riemannian structure on $\PP$ which pointwise corresponds to the Euclidean 
structure $\langle\cdot,\cdot\rangle$ 
on $T_x\PP$, $x\in \PP$ and  define three symplectic forms by 
$$\omega_{\II}(\cdot ,\cdot ) =  g(\cdot , \II\cdot ),\quad 
\omega_\JJ(\cdot ,\cdot ) =  g(\cdot , \JJ\cdot )\quad 
{\rm and} \quad \omega_\KK(\cdot ,\cdot ) =  g( \cdot , \KK\cdot ),$$ 
(in particular $\omega_\JJ \equiv \omega_\ominus$ and $\II\equiv J^T_{4n}$). 

Let $\Lambda$ be a $2n$-dimensional manifold and $\jmath :\Lambda\hookrightarrow \PP$
an embedding in the product manifold $\PP$. Consider a tubular neighborhood 
$\Lambda\subset U\subset \PP$ around $\Lambda$ being diffeomorphic to an 
open neighborhood around the zero section in $T^*\Lambda$ 
such that the projection $\pi:U\to \Lambda$ is well-defined and $\pi\circ\jmath = id_\Lambda$.

The following result characterizes the submanifolds $\Lambda$ which are adapted for 
constructing non-degenerated local symplectic maps. 
\begin{theorem}
\label{teo:main}
If the image $\Lambda\stackrel{\jmath}{\hookrightarrow} U\subseteq \PP$ is a
Lagrangian submanifold with respect to both $\omega_\II$ and $\omega_\JJ$ then: 
\begin{enumerate}
 \item it is a symplectic submanifold \footnote{ 
Note the similarity of the conditions on $\Lambda\subset \PP$ with those for
\emph{Special Lagrangian submanifolds} in K\"ahler or Calabi-Yau manifolds. 
See in particular \cite[Sec 8.1.1]{Joy07}.} with respect to $\omega_\KK$, 
  \item the kernel of the projection 
 $\pi:U\to \Lambda$  defines a local symplectic map by the equation 
 \begin{eqnarray}
   \pi_*(\JJ(v))=\pi_*(\II(v))=0 ,\qquad x\in\Lambda,\ v\in T_x\Lambda.
 \end{eqnarray}

\end{enumerate}

Such a map corresponds to the Cayley transformation of some 
Hamiltonian matrix $\bH\in End(TM)$.
 \end{theorem}
{\it Proof of 1.}
Let $x\in \Lambda\subset\PP$. For every tangent vector $v\in T_x\Lambda$
to the submanifold $\Lambda$, 
the vectors $\JJ(v)$ and $\II(v)$ belong to the normal bundle of $\Lambda$ in $\PP$
since it is Lagrangian 
for $\omega_\II$ and $\omega_\JJ$, 
\emph{i.e.} $\JJ(v),\II(v)\in (T_x\Lambda)^\perp$. 
In the same way, for every $u\in (T_x\Lambda)^\perp$ we have $\JJ(u),\II(u)\in T_x\Lambda$
and consequently 
\begin{eqnarray}
\II\circ\JJ(v) = - \JJ\circ\II(v) = \KK(v)\in T_x\Lambda. 
\end{eqnarray}
This shows that $T_x\Lambda$ and $(T_x\Lambda)^\perp$ are invariant under the action of $\KK$ which implies that
$\Lambda\subset \PP$ is a symplectic submanifold for $\omega_\KK$. 
Moreover, given the projection $\pi:U\subset\PP\to \Lambda$  we have $\pi_*(\JJ(v))=\pi_*(\II(v))=0$.
This is just the fact that $\ker \pi \equiv (T\Lambda)^\perp$.
$\hfill\square$

For proving the point {\it 2.} we need local coordinates and some additional elements. 
In fact, the proof is 
to explain how we can construct symplectic maps using Liouvillian forms. 
This is the subject of the following section.

\subsection{Liouvillian forms and symplectic maps}
Consider the same hypotheses of Theorem \ref{teo:main}.
For constructing symplectic maps using Liouvillian forms,
consider an element $v\in T_x\Lambda\subset T_x\PP$,
and search for primitive forms $\theta_\II$ and $\theta_\JJ$ such that $v\in\ker \theta_\JJ\cap\ker \theta_\II$. 
 Since $\Lambda$ is Lagrangian for $\omega_\II$ and $\omega_\JJ$, then 
 $\JJ(v),\II(v) \in (T_x\Lambda)^\perp$. Since $(T_x\Lambda)^\perp\equiv\ker \pi$, 
this implies $\pi_*\II(v) =\pi_*\JJ(v) = 0$. 
We will prove, in a contructive way, that the symplectic map is the solution of the equation 
$\pi_*\II(v)=0$. For free, we obtain that solving for a set of coordinates of one of the factor 
manifolds, gives the Cayley's transformation for some Hamiltonian matrix $\bH$.

Consider a local vector field $Z$ around $\Lambda$ being Liouville for both $\omega_\II$ and $\omega_\JJ$. By contraction,
we obtain the Liouvillian forms $\theta_\II=i_Z\omega_\II$ and 
$\theta_\JJ=i_Z\omega_\JJ$. 
Let $x\in \Lambda\subset\PP$ be a point and $v:=Z(x)\in T_x\Lambda$ 
the element of $Z$ on $T_x\Lambda$.
Then $v\in \ker \theta_\II\cap\ker \theta_\JJ$ by 
construction. We will construct a Liouville vector field $Z$
being suitable for constructing symplectic maps.
\begin{lemma}
  Let $\{x_i\}_{i=1}^{4n}$ be local coordinates on $\PP$.
  The ``expanding'' or ``Euler'' vector field $Z_0\in\Gamma(T\PP)$, given in these coordinates by 
  $Z_0=\fiz\sum_i x_i \frac{\partial}{\partial x_i}$, is 
  Liouville for all the three symplectic forms $\omega_\II$, $\omega_\JJ$ and $\omega_\KK$.
\end{lemma}
{\it Proof.} A direct verification shows that 
$$d\circ i_{(Z_0)} (\omega_C) = \omega_C, \qquad C\in\{\II,\JJ,\KK\},$$
and $\theta_C=i_{(Z_0)}\omega_C$ is a Liouvillian form for $(\PP,\omega_C)$.
$\hfill\square$

\begin{remark}
    The expanding vector field $Z_0$ is a degenerated case which corresponds to the identity map.
    In \cite{Jim15f} it is proved that the symplectic integrator contructed with 
    the expanding vector field corresponds to the mid point rule. 
\end{remark}

We proceed by looking for Liouville vector fields $Z\in \ker \theta_\II\cap \ker\theta_\JJ$ 
close to the expanding vector field. This is achieved by the addition of a component to 
the vector fields which is Hamiltonian with respect to $\omega_\II$ and $\omega_\JJ$. 
A vector field $Y=a_i(x)\frac{\partial}{\partial x_i}$ is Hamiltonian 
for $\omega_\II$ if $i_Y\omega_\II = - dF$ for a differentiable function 
$F:\PP\to\mathbb R$. In the same way, $Y$ is Hamiltonian for $\omega_\JJ$
if $i_Y\omega_\JJ = - dG$ for $G:\PP\to\mathbb R$. It means that 
$Y=\II\nabla F= \JJ\nabla G$ where $\nabla$ is the gradient associated to the 
Riemannian structure $g$ on $\PP$. Moreover, if $\psi:\PP\to\PP$ is a diffeomorphism 
taking $\II \stackrel{\psi}{\to} \JJ$ then $dF = \psi^*(dG)$.
However, this does not give information about the local structure of a 
vector field being Hamiltonian for both $\omega_\II$ and $\omega_\JJ$ at 
the same time.
For that we need to analyze the Jacobian matrix of $Y$ and without lost
of generality we consider that the components 
$a_i(x):\PP\to \mathbb R$ are linear
functions $a_i(x) = \sum_j A_{ij} x_j$.
The vector field $Y=\sum_i A_{ij}x_j\frac{\partial}{\partial x_i}$ is caracterized 
by the matrix $A=(A_{ij})$ and $Y$ is Hamiltonian for $\omega_\II$ if the matriz $A$ 
satisfies $A^T\II + \II A = 0$. Equivalently, $Y$ is Hamiltonian for $\omega_\JJ$ if
it satisfies $A^T\JJ + \JJ A = 0$.

\begin{lemma}
\label{lem:A}
   Let $S,R\in\mathbb M_{2n\times 2n}(\mathbb R)$ be a symmetric and a Hamiltonian matrix
   respectively, for the $2n$-dimensional symplectic manifold $(M,\omega)$. 
   Then the
   matrix $A\in\mathbb M_{4n\times 4n}(\mathbb R)$ given by
\begin{eqnarray}
A =   
    \left(\begin{smallmatrix}
        R & S\\
        - JSJ & - R^T
   \end{smallmatrix}\right).
   \label{eqn:A}
\end{eqnarray}
    is Hamiltonian for $(\PP,\omega_\II)$ and $(\PP,\omega_\JJ)$.
\end{lemma}
{\it Proof. } 
 The matrix $A$ is Hamiltonian for both $\omega_\II$ and $\omega_\JJ$ if it 
satisfies simultaneously: i) $A^T\II + \II A = 0$ and ii) $A^T\JJ + \JJ A = 0$.

Consider the matrix $ A = \left(\begin{smallmatrix}
        A_1 & A_2\\
        A_3 & A_4
   \end{smallmatrix}\right)$
and solving equation $A^T\II + \II A = 0$ gives the conditions 
$A_2=A_2^T$, $A_3=A_3^T$ and $A_4=-A_1^T$. It means $i)$ requires that $A_2$ and $A_3$ 
be symmetric and it relates $A_4$ with $A_1$.
On the other hand  the equation $A^T\JJ + \JJ A = 0$ gives the conditions 
$A_1^TJ + JA_1= 0$, $A_4^TJ+JA_4=0$ and $A_3 = -JA_2^TJ$. It means $ii)$ requires that $A_1$ and $A_4$
be Hamiltonian and it relates $A_2$ and $A_3$.
If we denote $R = A_1$ and $S= A_2$ then $A_3 = -JS^TJ$ and $A_4= -R^T$.
Finally, $R$ must be Hamiltonian for $\omega$ on $M$, and $S$ symmetric. This gives 
$A$ by expression (\ref{eqn:A}) which proves the lemma.
$\hfill\square$

If we consider that $A$ is not Hamiltonian for $\omega_\KK$ 
then $A^T\KK + \KK A\neq 0$. This produces the additional conditions 
$R\neq -R^T$ or $S\neq S^T$. Since $S=S^T$ is already a constraint from 
Lemma \ref{lem:A}, then $R$ cannot be antisymmetric. In particular, for 
$R$ a symmetric, Hamiltonian matrix for $(M,\omega)$ this problem has 
solutions.

In order to prove the second part of the main theorem, we need local 
coordinates for each one of the factors in the product manifold 
$\left(\{x_i\}_{i=1}^{2n},\{X_i\}_{i=1}^{2n}\right) \in M_1\times M_2 =:\PP$. In these coordinates,we have  
$$Z_0=\fiz\left( \begin{matrix} I & 0 \\ 0 & I \end{matrix} \right)\left( \begin{matrix} x \\ X \end{matrix} \right)\qquad{\rm and}\qquad Y=\fiz\left( \begin{matrix} R & S \\ -JSJ & -R^T \end{matrix} \right)\left( \begin{matrix} x \\ X \end{matrix} \right)$$
and the vector field $Z=Z_0+Y$ on $\PP$ is Liouville for $\omega_\II$ and $\omega_\JJ$. We add a $\fiz$ factor in $Y$ for simplify the computations.  The pointwise element $v=Z(x,X)$ is expressed in matricial form by 
\begin{eqnarray}
v=\fiz \left( \begin{smallmatrix} I + R & S \\ -JSJ & I-R^T\end{smallmatrix} \right)
\left( \begin{matrix} x \\ X \end{matrix} \right). 
\label{eqn:v}
\end{eqnarray}
We are in measure of proving second part of Theorem \ref{teo:main}.

{\it Proof of 2.}[Theorem \ref{teo:main}] 
Consider the vector $v=Z(x,X)$ given in matricial form by (\ref{eqn:v}).
This vector is tangent to $\Lambda$ by hypothesis and it is in the kernel of 
$\theta_\II$ and $\theta_\JJ$ by construction. 
Since $\Lambda\subset \PP$ is Lagrangian with respect to $\omega_\JJ$ then 
$\JJ^T(v)$ belongs to the normal bundle $(T_{(x,X)}\Lambda)^\perp$.
Applying the complex structure $\KK$ we have $(\KK\circ\JJ^T)(v)= \II(v)\in  (T_{(x,X)}\Lambda)^\perp$, with expression 

$$\II(v)=\fiz \left( \begin{matrix} - JSJ & I-R^T \\ -I -R & - S \end{matrix}\right) \left( 
 \begin{matrix} x \\ X \end{matrix} \right).$$
The equation $\pi_*\left( \II(v) \right) = 0$ in these local coordinates becomes
$$\left[-JSJ(x) + (I-R^T)(X) \right] + \left[ -(I+R)(x) - S(X)\right] = 0.$$
Rearranging  we obtain the matricial equation
$$\left[I - (R^T + S)\right]X = \left[I +(R+JSJ)\right]x.$$

Solving for $X$ is possible if $R^T + S$ is a non-exceptional matrix.\footnote{A matrix $A\in GL(n)$ is said 
to be \emph{non-exceptional} if $\det(I\pm A)\neq0$, where $I$ is the 
identity matrix in $GL(n)$.}
We consider the case where $R=R^T$ and $S=JSJ$, it means both matrices 
are symmetric and Hamiltonian.  Consequently, $\bH:=R^T+S = R + JSJ$ is well-defined
and it is a non-exceptional, Hamiltonian matrix for $(M,\omega)$. 
We solve for $X$ and we obtain 
$$X= (I-\bH)^{-1}(I+\bH)x.$$ 
The \emph{Cayley's transformation}  \cite{Wey46} assures that
the matrix $\bS=(I-\bH)^{-1}(I+ \bH)$ 
is symplectic if, and only if $\bH$ is Hamiltonian, and consequently the map 
$$x\mapsto (I-\bH)^{-1}(I+\bH)x$$ 
is a linear symplectic transformation. 
$\hfill\square$

Moreover, using the equation $\pi_*(\JJ(v))=0$, and following the same algebraic 
development, we arrive to the identity 
$X^\perp= (I-\bH)^{-1}(I+\bH)x^\perp$,
where $X^\perp=JX$ and $x^\perp=Jx$. 

The application in symplectic integration concerns the construction of
numerical integrators adapted to a given Hamiltonian system $(M,\omega,X_H)$. 
Denote by $\phi = (I-\tau\bH)^{-1}(I+\tau\bH)$ for small $\tau >0$ the symplectic map
obtained in Theorem \ref{teo:main}.
The implicit symplectic method \cite{KQ10,Jim15a}
\begin{eqnarray}
 z_\tau = z_0 + \tau X_H(\bar z)\quad{\rm where}\quad 
 \bar z=\fiz\left\{(z_0 + z_\tau) +\tau \bH (z_h - z_0) \right\},
\end{eqnarray}
integrates the Hamiltonian vector field 
$\dot\zeta = X_{\hat H}(\zeta)$ for an alternative Hamiltonian function 
$\hat H:M\to\mathbb R$, known as the ``surrounding Hamiltonian'' 
in backward error analysis \cite{Rei99}.
The surrounding Hamiltonian is related to the original Hamiltonian by 
$\hat H = H\circ\phi^{-1}$.

\end{document}